\documentclass{amsart}
\usepackage{amsfonts}
\usepackage{amsmath}
\usepackage{amssymb}

\theoremstyle{plain}

\newtheorem{remark}{Remark}

\numberwithin{equation}{section}

\input{tcilatex}

\begin{document}
\title[Zeros of the Riemann zeta function]{\textbf{Conditional and
Unconditional Large Gaps Between the Zeros of the Riemann Zeta-Function}}
\author{\textbf{S. H. Saker}}
\address{{\small Department of Mathematics Skills, PYD, King Saud
University, Riyadh 11451, } {\small Saudi Arabia, Department of Math.,
Faculty of Science, Mansoura University, Mansoura 35516, Egypt.}}
\email{{\small shsaker@mans.edu.eg, mathcoo@py.ksu.edu.sa}}
\maketitle

\begin{abstract}
In this paper, first by employing inequalities\ derived from the Opial
inequality due to David Boyd with best constant, we will establish new
unconditional lower bounds for the gaps between the zeros of the Riemann
zeta function. Second, on the hypothesis that the moments of the Hardy $Z-$%
function and its derivatives are correctly predicted, we establish some
explicit formulae for the lower bounds of the gaps between the zeros and use
them to establish some new conditional bounds. In particular it is proved
that the consecutive nontrivial zeros often differ by at least $6.1392$
(conditionally) times the average spacing. This value improves the value $%
4.71474396$ that has been derived in the literature.
\end{abstract}

\section{Introduction}

The Riemann zeta function $\zeta (s)$ is defined on $\{s\in \mathbb{C}:\func{%
Re}(s)>1\}$ by the series\ 
\begin{equation*}
\zeta (s):=1+\frac{1}{2^{s}}+\frac{1}{3^{s}}+\frac{1}{4^{s}}+\frac{1}{5^{s}}%
+...,
\end{equation*}%
which converges in the region described by the Cauchy integral test. It is
of fundamental importance because it can also be represented just in terms
of the primes. This representation is given by 
\begin{equation*}
\zeta (s):=\dprod_{p}\left( 1-\frac{1}{p^{s}}\right) ^{-1},\ \ \ \text{for\ }%
\ \ \ \func{Re}s>1,
\end{equation*}%
where the product is taken over all prime numbers. Thus its analytic
properties are related to the distribution of prime numbers. Among the
integers the primes appear to be scattered at random. It is known that they
are infinite in number, but there is no useful formula which generate them.
However, on average they obey simple laws. For example, the prime number
theorem states that the number of primes which occur up to a given integer $%
X $, $\pi (X)$, is approximately $X/\log (X)$, the approximation getting
better as $X$ increases. The actual numbers found for different $X$ will
fluctuate about this value.

\bigskip

Riemann gave an exact formula for the counting function $\pi (X)$, in which
fluctuations about the average are related to the values of $s$ for which $%
\zeta (s)=0.$ These are isolated points in the complex plan. In \cite{jamal}
the authors presented a numerical study of Riemann's formula for the
oscillating part of the density of the primes and their integer powers. The
formula consists of an infinite series of oscillatory terms, one for each
zero of the zeta function on the critical line, and was derived by Riemann
in his paper on primes, assuming the Riemann hypothesis. They also showed
that high-resolution spectral lines can be generated by the truncated series
at all integer powers of primes and demonstrate explicitly that the relative
line intensities are correct. They then derived a Gaussian sum rule for
Riemann's formula and used to analyze the numerical convergence of the
truncated series.

Riemann conjectured that all nontrivial (non-real) zeros are distributed
symmetrically with respect to the critical line $\func{Re}s=1/2$ and the
real axis. This is the Riemann hypothesis. Riemann showed that the
zeta-function satisfies a functional equation of the form%
\begin{equation}
\pi ^{-s/2}\Gamma (\frac{s}{2})\zeta (s)=\pi ^{-(1-s)/2}\Gamma (\frac{1-s}{2}%
)\zeta (1-s),  \label{F}
\end{equation}%
where $\Gamma $ is the Euler gamma function. Clearly, there are no zeros in
the half-plane of convergence $\func{Re}(s)>1,$ and it is also known that $%
\zeta (s)$ does not vanish on the line $\func{Re}(s)=1.$ In the
negative-half plane $\zeta (s)$ and its derivative are oscillatory and from
the functional equation there exist-so called trivial (real) zeros at $s=-2m$
for any positive integer $m$ (corresponding to the poles of the appearing
Gamma-factors). It is conjectured that all or at least almost all nontrivial
zeros of the zeta-function are simple, see \cite{CG} and \cite{C2}.

\bigskip

Reimann's connection between the nontrivial zeros and the primes has
particularly interesting form: it bears a striking resemblance to the
Gutzwiller formula, with the zeros behaving like energy levels and the
primes labelling the periodic orbits of some chaotic classical system.
Montgomery \cite{M1} studied the distribution of pairs of nontrivial zeros $%
1/2+i\gamma $ and $1/2+i\gamma ^{^{\prime }}$ and conjectured, for fixed $%
\alpha ,$ $\beta \,$\ satisfying $0<\alpha <\beta ,\,$that 
\begin{eqnarray*}
&&\lim_{T\rightarrow \infty }\frac{1}{N(T)}\#\left\{ 0<\gamma ,\gamma
^{^{\prime }}<T:\alpha \leq \frac{\gamma ^{^{\prime }}-\gamma ^{^{\prime }}}{%
(2\pi /\log T)}\leq \beta \right\} \\
&=&\int_{\alpha }^{\beta }\left( 1-\left( \frac{\sin \pi x}{\pi x}\right)
^{2}\right) dx.
\end{eqnarray*}%
This so-called pair correlation conjecture plays a complementary role to the
Riemann hypothesis. This conjecture implies the essential simplicity
hypothesis that almost all zeros of the zeta-function are simple. On the
other hand the integral on the right hand side is the same as the one
observed in the two point correlation of the eigenvalues which are the
energy levels of the corresponding Hamiltonian that are usually not known
with uncertainty. This observation is due to Dyson and it restored some hope
in an old idea of Hilbert and Polya that the Riemann hypothesis follows from
the existence of a self-adjoint Hermitian operator whose spectrum of
eigenvalues correspond to the set of nontrivial zeros of the zeta function.

\bigskip

Odlyszko \cite{Odl} published the results of a remarkable series of computer
calculations of the zeros which showed that they were the same as those of
large Hermitian matrices with randomly picked entries. These suggests that
the zeros might well be the energy levels of some as yet unidentified
quantum system whose classical motion is chaotic, and without symmetry under
time reversal. The connections to quantum chaos and semiclassical physics
are discussed in \cite{jamal}. So that the distribution of zeros of the
Riemann zeta-function is of fundamental importance in number theory as well
as in physics.

\bigskip

The number $N(t)$ of the non-trivial zeros of $\zeta (s)$ with ordinate in
the interval $[0,$ $T]$ is asymptotically given by the Riemann-von Mangoldt
formula (see \cite{RS}) 
\begin{equation*}
N(T)=\frac{T}{2\pi }\log (\frac{T}{2\pi e})+O(\log T).
\end{equation*}%
Consequently there are infinitely many nontrivial zeros, all of them lying
in the critical strip $0<\func{Re}s<1,$ and the frequency of their
appearance is increasing as $T\rightarrow \infty .$ Assume that $(\beta
_{n}+i\gamma _{n})$ are the zeros of $\zeta (s)$ in the upper half-plane
(arranged in non-decreasing order and counted according multiplicity) and $%
\gamma _{n}\leq \gamma _{n+1}$ are consecutive ordinates of all zeros.
Define 
\begin{equation}
\lambda :=\lim \sup_{n\rightarrow \infty }\frac{(\gamma _{n+1}-\gamma _{n})}{%
(2\pi /\log \gamma _{n})},\text{ and }\mu :=\lim \inf_{n\rightarrow \infty }%
\frac{(\gamma _{n+1}-\gamma _{n})}{(2\pi /\log \gamma _{n})},  \label{Q}
\end{equation}%
where $\left( 2\pi /\log \gamma _{n}\right) $ is the average spacing between
zeros. The values of $\lambda $ and $\mu $ have received a great deal of
attention. In fact, important results have been obtained by some authors. It
generally conjectured that 
\begin{equation}
\mu =0,\text{ \ \ and \ \ \ }\lambda =\infty .  \label{b}
\end{equation}%
As mentioned by Montogomery \cite{M1} it would be interesting to see how
numerical evidence compare with the above conjectures. Now, several results
has been obtained, however the failure of Gram's low (see \cite{H})
indicates that the asymptotic behavior is approached very slowly. Thus the
numerical evidence may not be particularly illuminating. So that any
numerical values of $\mu $ and $\lambda $ may be help in proving (\ref{b}),
which is one of our aims in this paper. Selberg \cite{Selberg} proved that $%
0<\mu <1<\lambda $ and the average of $r_{n}$ is $1.$ Mueller \cite{33}
obtained $\lambda >1.9$ assuming the Riemann hypothesis. Montogomery and
Odlyzko \cite{Montgomery} showed, assuming the Riemann hypothesis, that $%
\lambda >1.9799$\ and\ $\mu <0.5179.$ Conrey, Ghosh and Gonek \cite{Conrey1}
showed that, if the Riemann hypothesis is true, then $\lambda >2.337,$ and\ $%
\mu <0.5172.$ Bui, Milinovich and Ng \cite{BMN} obtained $\lambda >2.69,$
and $\mu <0.5155$ assuming the Riemann hypothesis. Conrey, Ghosh and Gonek 
\cite{Conrey2} obtained a new lower bound and proved that $\lambda >2.68$
assuming the generalized Riemann hypothesis for the zeros of the Dirichlet $%
L-$ functions. Ng in \cite{Ng} proved that $\lambda >3$ assuming the
generalized Riemann hypothesis for the zeros of the Dirichlet $L-$functions.
Bui\cite{Bui} proved that $\lambda >3.0155$ assuming the generalized Riemann
hypothesis for the zeros of the Dirichlet $L-$functions. The main results in 
\cite{Bui, BMN, Conrey2, Ng} are based on the idea of Mueller \cite{33}.

\bigskip

Hall \cite{Hall1} supposed that the sequence of distinct positive zeros of
the Riemann zeta-function $\zeta (\frac{1}{2}+it)$ which arranged in
non-decreasing order and counted according multiplicity is given by $%
\{t_{n}\}$ and defined%
\begin{equation}
\Lambda :=\lim \sup_{n\rightarrow \infty }\frac{t_{n+1}-t_{n}}{(2\pi /\log
t_{n})}\text{,}  \label{R2}
\end{equation}%
which is the quantity in (\ref{Q}) where only zeros $\frac{1}{2}+it_{n}$ on
the critical line with the idea that this could be bounded from below
unconditionally. Note that the Riemann hypothesis implies that the $t_{n}$
corresponded to the positive ordinates of non-trivial zeros of the zeta
function, i.e., $N(T)\sim \left( T\log T\right) /2\pi .$ The average spacing
between consecutive zeros with ordinates of order $T$ is $2\pi /\log (T)$
which tends to zero as $T\rightarrow \infty .$

\bigskip

Hall \cite{Hall3} showed that $\Lambda \geq \lambda $, and the lower bound
for $\Lambda $ bear direct comparison with such bounds for $\lambda $
dependent on the Riemann hypothesis, since if this were true the distinction
between $\Lambda $ and $\lambda $ would be nugatory. Of course $\Lambda \geq
\lambda $ and the equality holds if the Riemann hypothesis is true. So that
if the Riemann hypothesis is true, we see that any improvement of $\Lambda $
(unconditionally) will lead to the improvement of $\lambda $ and vice versa.
The behavior of $\zeta (s)$ on the critical line is reflected by the Hardy $%
Z-$function $Z(t)$ as a function of a real variable, defined by 
\begin{equation}
Z(t)=e^{i\theta (t)}\zeta (\frac{1}{2}+it),\text{ where }\theta (t):=\pi
^{-it/2}\frac{\Gamma (\frac{1}{4}+\frac{1}{2}it)}{\left\vert \Gamma (\frac{1%
}{4}+\frac{1}{2}it)\right\vert }.  \label{Z}
\end{equation}%
Also it follows that $Z(t)$ is an infinitely often differentiable and real
for real $t$ and moreover $\left\vert Z(t)\right\vert =\left\vert \zeta
(1/2+it)\right\vert $. Consequently, the zeros of $Z(t)$ correspond to the
zeros of the Riemann zeta-function on the critical line. In \cite{Hall2}
Hall proved a Wirtinger-type inequality and used the moment 
\begin{equation}
\int_{0}^{T}Z^{4}(t)dt=\frac{1}{2\pi ^{2}}T\log ^{4}(T)+O(T\log ^{3}T)\text{%
, }  \label{hall2}
\end{equation}%
due to Ingham \cite{Ingham}) and the moments%
\begin{equation}
\int_{0}^{T}(Z^{^{\prime }}(t))^{4}dt=\frac{1}{1120\pi ^{2}}T\log
^{8}(T)+O(T\log ^{7}T),  \label{hall3}
\end{equation}%
\begin{equation}
\int_{0}^{T}Z^{2}(t)(Z^{^{\prime }}(t))^{2}dt=\frac{1}{120\pi ^{2}}T\log
^{6}(T)+O(T\log ^{5}T),  \label{hall4}
\end{equation}%
due to Conrey \cite{Conrey3}, and obtained unconditionally that $\Lambda
\geq 2.3452.$

\bigskip

The moments $I_{k}(T)$ of the Hardy $Z-$function $Z(t)$ and the moments $%
M_{k}(T)$ of its derivative are defined by 
\begin{equation*}
I_{k}(T):=\int_{0}^{T}\left\vert Z(t)\right\vert ^{2k}dt,\text{ and }%
M_{k}(T):=\int_{0}^{T}\left\vert Z^{^{\prime }}(t)\right\vert ^{2k}dt.
\end{equation*}%
For positive real numbers $k$, it is believed that $I_{k}(T)\sim C(k)$ $%
T\left( \log T\right) ^{k^{2}}$and $M_{k}(T)\sim L(k)T\left( \log T\right)
^{k^{2}+2k}$ for positive constants $C_{k}$ and $L_{k}$ will be defined
later. Keating and Snaith \cite{KS} based on considerations from random
matrix theory conjectured that%
\begin{equation}
I_{k}(T)\sim a(k)b(k)T\left( \log T\right) ^{k^{2}},  \label{A1}
\end{equation}%
where $a(k)$ is a product over the primes which is defined by%
\begin{equation*}
a(k):=\tprod_{p}(\left( 1-\frac{1}{p^{2}}\right) \dsum_{m=0}^{\infty }\left( 
\frac{\Gamma (m+k)}{m!\Gamma (k)}\right) ^{2}p^{-m},
\end{equation*}%
and 
\begin{equation}
b(k):=\dprod_{j=0}^{k-1}\frac{j!}{(j+k)!}.  \label{B}
\end{equation}%
Using the relation (\ref{B}) one can obtain the value of $b(k)$ for any real
positive number $k.$ Conrey, Rubinstein and Snaith \cite{C} conjectured that 
\begin{equation}
M_{k}(T)\sim a(k)c(k)T\left( \log T\right) ^{k^{2}+2k},  \label{A2}
\end{equation}%
where 
\begin{equation*}
c(k):=(-1)^{\frac{k(k+1)}{2}}\dsum_{m\in P_{O}^{k+1}(2k)}^{k}\left( 
\begin{array}{c}
2k \\ 
m%
\end{array}%
\right) \left( \frac{-1}{2}\right) ^{m_{0}}\left( \dprod_{i=1}^{k}\tfrac{1}{%
(2k-i+m_{i})!}\right) M_{i,j},
\end{equation*}%
where 
\begin{equation}
M_{i,j}:=\left( \dprod_{1\leq i,j\leq k}^{k}(m_{j}-m_{i}+i-j)\right) ,
\label{M}
\end{equation}%
and $P_{O}^{k+1}(2k)$ denotes the set of partitions $m=(m_{0},...,m_{k})$ of 
$2k$ into nonnegative parts. In this paper, we will determine the values of $%
b(k)/c(k)$ for $k=1,2,...,15$ that we will use derive the new conditional
lower bounds for $\Lambda .$

\bigskip

Hall in \cite{Hall3, Hall8} used the moments of mixed powers of the form 
\begin{equation}
\int_{0}^{T}Z^{2k-2h}(t)(Z^{^{\prime }}(t))^{2h}dt\sim C(h,k)T\left( \log
T\right) ^{k^{2}+2h},  \label{M1}
\end{equation}%
where $Z(t)$ is the Hardy $Z-$function, and $k\in N$, $0\leq h\leq k$ and a
complicated variation problem together with a new Wirtinger-type inequality
designed exclusively for this problem and obtained some conditional lower
bounds of $\Lambda $. The moments in (\ref{M1}) has been predicted by Random
Matrix Theory (RMT) by Hughes \cite{Hughes} who stated an interesting
conjecture on the moments of the zeta function and its derivatives at its
zeros subject to the truth of Riemann's hypothesis when the zeros are
simple. This conjecture includes for fixed $k>-3/2$ the asymptotes formula
of the moments of the higher order of the Riemann zeta function and its
derivative. We suppose further that if $k$ is a fixed positive integer and $%
h\in \lbrack 0,$ $k]$ is an integer then the formula 
\begin{equation}
\int_{0}^{T}Z^{2k-2h}(t)(Z^{^{\prime }}(t))^{2h}dt\sim a(k)b(h,k)T\left(
\log T\right) ^{k^{2}+2h},  \label{Hu}
\end{equation}%
holds. Note that this was predicted by Keating and Snaith \cite{KS} in the
case when $h=0$, with wider range $\func{Re}(k)>-1/2$ and by Hughes \cite%
{Hughes} in the range $\min (h,k-h)>-1/2,$ $a(k)$ is a product over the
primes and $b(h,k)$ is rational: indeed for integral $h$, it is obtained
that 
\begin{equation}
b(h,k)=b(k)\left( \frac{\left( 2h\right) !}{8^{h}h!}\right) H(h,k),
\label{BG}
\end{equation}%
where $H(h,k)$ is an explicit rational function of $k$ for each fixed $h.$
The functions $H(h,k)$ as introduced by Hughes \cite{Hughes} are given in
the following table where $K=2k$:

\begin{center}
\begin{equation*}
\begin{tabular}[t]{|l|}
\hline
$H(0,k)=1,$ $H(1,k)=\frac{1}{K^{2}-1},$ $H(2,k)=\frac{1}{(K^{2}-1)(K^{2}-9)}$
\\ \hline
$H(3,k)=\frac{1}{(K^{2}-1)^{2}(K^{2}-25)},$ $H(4,k)=\frac{K^{2}-33}{%
(K^{2}-1)^{2}(K^{2}-9)(K^{2}-25)(K^{2}-49)}$ \\ \hline
$H(5,k)=\frac{K^{4}-90K^{2}+1497}{%
(K^{2}-1)^{2}(K^{2}-9)^{2}(K^{2}-25)(K^{2}-49)(K^{2}-81)},$ \\ \hline
$H(6,k)=\frac{K^{6}-171K^{4}+6867K^{2}-27177}{%
(K^{2}-1)^{3}(K^{2}-9)^{2}(K^{2}-25)(K^{2}-49)(K^{2}-81)(K^{2}-121)},$ \\ 
\hline
$H(7,k)=\frac{K^{8}-316K^{6}+30702K^{4}-982572K^{2}+6973305}{%
(K^{2}-1)^{3}(K^{2}-9)^{2}(K^{2}-25)^{2}(K^{2}-49)(K^{2}-81)(K^{2}-121)(K^{2}-169)%
},$ \\ \hline
\end{tabular}%
\end{equation*}%
\textit{Table 1. The values of }$H(h,k)$\textit{, where }$K=2k.$
\end{center}

This sequence continuous, and it is believed that both the nominator and
denominator are polynomials in $k^{2}$, moreover that the denominator is
actually (see \cite{D}) 
\begin{equation}
\dprod_{a\text{ }odd>0}\left\{ (K^{2}-a^{2})^{\alpha (a,h)}:\alpha (a,h)=%
\frac{4h}{a+\sqrt{a^{2}+8h}}\right\} .  \label{Monic}
\end{equation}%
Using the equation (\ref{BG}) and the definitions of the functions $H(h,k)$,
we can obtain the values of $b(0,k)/b(k,k)$ for $k=1,2,...,7.$ As indicated
in \cite{Hall8} Hughes \cite{Hughes} conjectured the first four functions
and then writes that numerical experiment suggests the next three. Hall \cite%
{Hall8} shown that in the case when $h=3$, ($H(3,k))$ requires adjustment to
fit with (\ref{Monic}) in that extra factor $K^{2}-9$ should be introduced
in both the nominator and denominator. To use (\ref{Hu}) Hall \cite{Hall3}
proved a new generalized Wirtinger-type inequality of the form 
\begin{equation}
\int_{0}^{\pi }H\left( y^{^{\prime }}(t)/y(t)\right) y^{2k}(t)dt\geq
(2k-1)L\int_{0}^{\pi }y^{2k}(t)dt,  \label{W2}
\end{equation}%
where $y(t)\in C^{2}[0,$ $\pi ],$ $y(0)=y(\pi )=0$, $L=L(k,H)$ is determined
from the solution of the equation 
\begin{equation*}
\int_{0}^{\infty }\frac{G^{^{\prime }}(u)}{G(u)+(2k-1)L}\frac{du}{u}=k\pi 
\text{, for }k\in \mathbb{N},
\end{equation*}%
where $G(u):=uH^{^{\prime }}(u)-H(u)$, $H(u)$ be an even function,
increasing, strictly convex on $\mathbb{R}^{+}$ and satisfies $%
H(0)=H^{^{\prime }}(0)=0,$and $uH^{^{\prime \prime }}(u)\rightarrow 0$ as $%
u\rightarrow 0.$ The inequality (\ref{W2}) is proved by using the calculus
of variation which depends on the minimization of the integral on the left
hand side subject to the constrains $y(0)=0$ and $\int_{0}^{\pi
}y^{2k}(t)dt=1.$ Assuming that\ (\ref{Hu}) is correctly predicted, Hall
employed the inequality (\ref{W2}) when 
\begin{equation*}
H(u):=\dsum_{h=1}^{k}\frac{2k-1}{2h-1}\left( 
\begin{array}{c}
h \\ 
k%
\end{array}%
\right) \upsilon _{h}u^{2h}\text{, \ \ }\upsilon _{h}\geq 0\text{, \ \ }%
\upsilon _{k}=1,
\end{equation*}%
and obtained an explicit formula for $\Lambda (k)$ which is given by 
\begin{equation*}
\Lambda ^{2}\geq X,
\end{equation*}%
where $X$ is the real positive root of the equation 
\begin{equation*}
\sum_{h=1}^{k}\frac{2k-1}{2h-1}\left( 
\begin{array}{c}
k \\ 
h%
\end{array}%
\right) R(h,k)\upsilon _{h}X^{h}-(2k-1)\lambda (\upsilon _{1},\upsilon
_{1},...,\upsilon _{k-1})R(0,k)=0,
\end{equation*}%
and 
\begin{equation*}
R(h,k):=4^{k-h}\frac{b(h,k)}{b(k,k)}.
\end{equation*}%
He then derived a new value of $\Lambda $ (when $k=3)$ which is given by 
\begin{equation}
\Lambda \geq \sqrt{7533/901}=2.8915.  \label{D2}
\end{equation}%
The main challenge in \cite{Hall3} was to maximize $X=\kappa ^{2}$ (which is
not an easy task) where $X$ satisfies the equation 
\begin{equation*}
27X^{3}+385\mu X^{2}+10395\vartheta X-121275L=0,
\end{equation*}%
and $L$ obtained form the equation 
\begin{equation*}
\int_{-\infty }^{\infty }\frac{x^{4}+2\mu x^{2}+\upsilon }{x^{6}+3\mu
x^{4}+3\upsilon x^{2}+L}dx=\pi .
\end{equation*}%
Hall \cite{Hall4} simplified the calculations in \cite{Hall3} and converted
the problem into one of the classical theory of equations involving
Jacobi-Schur functions and proved that $\Lambda (4)\geq 3.392272,$ $\Lambda
(5)\geq 3.858851,$ and $\Lambda (6)\geq 4.2981467.$ Hall \cite{Hall8}
developed the theory set used in \cite{Hall4} and proved that $\Lambda
(7)\geq 4.215007$ assuming that (\ref{Hu}) is correctly predicted. The
improvement of this value as obtained in \cite{Hall8} is given by $\Lambda
(7)\geq 4.71474396$ assuming that (\ref{Hu}) is correctly predicted. The
question now is: If it is possible to employ new inequalities with best
constants to find new explicit formulae for the gaps and use them to find
new series of the lower bounds?

\bigskip

The paper gives an affirmative answer to this question. In fact, we will
derive new unconditional and conditional lower bounds for $\Lambda .$ The
main results will be proved by employing two inequalities derived from the
Opial inequality with a best constant due to David Boyd \cite{Boyd} who
applied a variational technique to reduce the determination of the best
constant to a nonlinear eigenvalue problem for an integral operator.

\bigskip

The main results will be proved in the next section which is organized as
follows: First, we derive new unconditional lower bounds for $\Lambda $.
Second on the hypothesis that the moments of the Hardy $Z-$function and its
derivatives are correctly predicted, we establish new explicit formulae of
the gaps between the zeros and establish some lower bounds for $\Lambda $.
In particular, we will prove that $\Lambda \geq 6.1392$ which improves the
value $\Lambda \geq 4.71474396.$

\section{Main Results}

Before we state and prove the main results, we derive some inequalities from
the Opial inequality due to David Boyd \cite{Boyd} that we will use in this
section. The Opial inequality due to David Boyd \cite{Boyd} is presented in
the following theorem.

\bigskip

\textbf{Theorem A.} \textit{If }$y\in C^{1}[a,$\textit{\ }$b]$\textit{\ with 
}$y(a)=0$\textit{\ (or }$y(b)=0)$\textit{, then}%
\begin{equation}
\int_{a}^{b}\left\vert y(t)\right\vert ^{p}\left\vert y^{^{\prime
}}(t)\right\vert ^{q}dt\leq K(p,q,r)(b-a)^{r-q}\left( \int_{a}^{b}\left\vert
y^{^{\prime }}(t)\right\vert ^{r}dt\right) ^{\frac{p+q}{r}},  \label{Boyd}
\end{equation}%
\textit{where }$p>0$\textit{, }$r>1$\textit{, }$0\leq q<r$\textit{,} 
\begin{eqnarray*}
K(p,q,r) &:&=\frac{\left( r-q\right) p^{p}}{(r-1)(p+q)}\beta ^{p+q-r}\left(
I(p,q,r)\right) ^{-p},\text{ } \\
\beta &:&=\left\{ \frac{p(r-1)+(r-q)}{(r-1)(p+q)}\right\} ^{\frac{1}{r}},
\end{eqnarray*}%
\textit{and}%
\begin{equation*}
I(p,q,r):=\int_{0}^{1}\left\{ 1+\frac{r(q-1)}{r-q}t\right\}
^{-(p+q+rp)/rp}[1+(q-1)t]t^{1/p-1}dt.
\end{equation*}%
First, we will derive a new inequality from the Opial inequality (\ref{Boyd}%
) of the from (\ref{W2}) which allows us to use the moments (\ref{hall3})
and (\ref{hall4}) to derive the new unconditional lower bound of $\Lambda .$
For a special case of Theorem A, when $r=p+q$, we have 
\begin{equation}
\int_{a}^{b}\left\vert y(t)\right\vert ^{p}\left\vert y^{^{\prime
}}(t)\right\vert ^{q}dt\leq K(p,q,p+q)(b-a)^{p}\left( \int_{a}^{b}\left\vert
y^{^{\prime }}(t)\right\vert ^{p+q}dt\right) ,  \label{KV}
\end{equation}%
where%
\begin{equation}
K(p,q,p+q):=\frac{q(p+q)^{p-1}}{\left( pL(p,q)+q\right) ^{p}},\text{ }q\neq 0%
\text{,}  \label{KV1}
\end{equation}%
and 
\begin{equation}
L(p,q):=\int_{0}^{1}\left( \frac{1}{1-\lambda s^{p}}\right) ds\text{, \
where\ }\lambda =\frac{(p+q)(q-1)}{(p+q-1)q}.  \label{Kv2}
\end{equation}%
The inequality (\ref{KV}) has immediate application to the case where $%
y(a)=y(b)=0$. Choose $c=(a+b)/2$ and apply (\ref{B1}) to $[a,c]$ and $[c,b]$
and then add to obtain%
\begin{eqnarray*}
&&\int_{a}^{b}\left\vert y(t)\right\vert ^{p}\left\vert y^{^{\prime
}}(t)\right\vert ^{q}dt \\
&\leq &K(p,q,p+q)(\frac{b-a}{2})^{p}\left( \int_{a}^{c}\left\vert
y^{^{\prime }}(t)\right\vert ^{p+q}dt+\int_{c}^{b}\left\vert y^{^{\prime
}}(t)\right\vert ^{p+q}dt\right) \\
&\leq &K(p,q,p+q)(\frac{b-a}{2})^{p}\left( \int_{a}^{b}\left\vert
y^{^{\prime }}(t)\right\vert ^{p+q}dt\right) .
\end{eqnarray*}%
So that if $y(0)=y(\pi )=0$, we have%
\begin{equation}
\int_{0}^{\pi }\left\vert y(t)\right\vert ^{p}\left\vert y^{^{\prime
}}(t)\right\vert ^{q}dt\leq K(p,q,p+q)(\frac{\pi }{2})^{p}\left(
\int_{0}^{\pi }\left\vert y^{^{\prime }}(t)\right\vert ^{p+q}dt\right) .
\label{Kv3}
\end{equation}%
If we choose $p=2$ and $q=2$, we get that 
\begin{equation}
\int_{0}^{\pi }y^{2}(t)(y^{^{\prime }}(t))^{2}dt\leq K(2,2,4)(\frac{\pi }{2}%
)^{2}\left( \int_{0}^{\pi }\left( y^{^{\prime }}(t)\right) ^{4}dt\right) .
\label{Kv4}
\end{equation}%
Using the definition of $K$, we see that $K(2,2,4)=0.346\,13,$ where we used
the value of 
\begin{equation*}
L(2,2)=\int_{0}^{1}\left( \frac{1}{1-\frac{2}{3}s^{2}}\right) ds=1.4038.
\end{equation*}%
So that the inequality (\ref{Kv4}) becomes 
\begin{equation}
\left( \int_{0}^{\pi }\left( y^{^{\prime }}(t)\right) ^{4}dt\right) \geq 
\frac{4}{\left( 0.346\,13\right) \pi ^{2}}\int_{0}^{\pi
}y^{2}(t)(y^{^{\prime }}(t))^{2}dt.  \label{Kv5}
\end{equation}%
By a suitable linear transformation, we deduce that if $y\in C^{1}[a,b]$%
\textit{\ }with $y(a)=0=y(b),$ then we have%
\begin{equation}
\int_{a}^{b}\left( \frac{b-a}{\pi }\right) ^{2}\left( y^{^{\prime
}}(t)\right) ^{4}dt\geq \frac{4}{\left( 0.346\,13\right) \pi ^{2}}%
\int_{a}^{b}y^{2}(t)(y^{^{\prime }}(t))^{2}dt.  \label{Kv6}
\end{equation}%
In the following, assuming the Riemann hypothesis, we will apply the
inequality (\ref{Kv6}) and using the moments (\ref{hall3}) and (\ref{hall4})
to find a new unconditional value of $\Lambda $. One can see that the value
that we will establish\ does not improve the obtained values, but the
technique is a simple one and depends only on the application of an
inequality derived from the well-known Opial inequality. Note that, as
mentioned by Hall, when the Riemann hypothesis is true we have $\lambda
=\Lambda .$

\bigskip

\textbf{Theorem 2.1}.\textit{\ Let }$\varepsilon (T)\rightarrow 0$\textit{\
in such a way that }$\varepsilon (T)\log T\rightarrow \infty .$\textit{\
Then for sufficiently large }$T$\textit{, there exists an interval contained
in }$[T,(1+\varepsilon (T))T]$\textit{\ which is free of zeros of }$Z(t)$%
\textit{\ and having length at least }%
\begin{equation*}
\sqrt{\frac{112}{12\left( 0.34613\right) \pi ^{2}}}\left\{ 1+O\left( \frac{1%
}{\varepsilon (T)\log T}\right) \right\} \frac{2\pi }{\log T}.
\end{equation*}%
\textit{Thus}%
\begin{equation}
\Lambda \geq \frac{1}{\pi }\sqrt{\frac{112}{12\left( 0.34613\right) }}=1.6529%
\text{.}  \label{h1}
\end{equation}%
\textbf{Proof. }We follow the arguments in \cite{Hall2} to prove our
theorem. Suppose that $t_{l}$ is the first zero of $Z(t)$ not less than $T$
and $t_{m}$ the last zero not greater than $(1+\varepsilon )T$ where $%
\varepsilon (T)\rightarrow 0$ in such a way that $\varepsilon (T)\log
T\rightarrow \infty .$ Suppose further that for $l\leq n<m$, we have 
\begin{equation}
L_{n}=t_{n+1}-t_{n}\leq \frac{2\pi \kappa }{\log T}.  \label{kv3}
\end{equation}%
Applying the inequality (\ref{Kv6}) with $y(t)=Z(t)$, we have 
\begin{equation*}
\int_{t_{n}}^{t_{n+1}}\left[ \left( \frac{L_{n}}{\pi }\right) ^{4}\left(
Z^{^{\prime }}(t)\right) ^{4}-\frac{4}{\left( 0.346\,13\right) \pi ^{2}}%
Z^{2}(t)(Z^{^{\prime }}(t))^{2}\right] dt\geq 0.
\end{equation*}%
Since the inequality remains true if we replace $L_{n}/\pi $ by $2\kappa
/\log T$, we have 
\begin{equation}
\int_{t_{n}}^{t_{n+1}}\left[ \left( \frac{2\kappa }{\log T}\right)
^{4}\left( Z^{^{\prime }}(t)\right) ^{4}-\frac{4}{\left( 0.346\,13\right)
\pi ^{2}}Z^{2}(t)(Z^{^{\prime }}(t))^{2}\right] dt\geq 0.  \label{kv4}
\end{equation}%
Summing (\ref{kv4}) over $n,$ using (\ref{hall3}) and (\ref{hall4}), we
obtain 
\begin{eqnarray*}
&&\frac{1}{1120\pi ^{2}}\left( \frac{2\kappa }{\log T}\right) ^{2}T\log
^{8}(T)+O(T\log ^{7}) \\
&&-\frac{4}{\left( 0.346\,13\right) \pi ^{2}}\frac{1}{120\pi ^{2}}T\log
^{6}(T)+O(T\log ^{5}T) \\
&=&\frac{\left( 2\kappa \right) ^{2}}{1120\pi ^{2}}T\log ^{6}(T)+O(T\log
^{7}T) \\
&&-\frac{4}{\left( 0.346\,13\right) \pi ^{2}}\frac{1}{120\pi ^{2}}(T\log
^{6}T)+O(T\log ^{5}T).
\end{eqnarray*}%
Follows the proof of Theorem 1 in \cite{Hall2}, we obtain 
\begin{equation*}
\kappa ^{2}\geq \frac{112}{12\left( 0.346\,13\right) \pi ^{2}}%
+O(1/\varepsilon (T)\log T).
\end{equation*}%
Then, we have (noting $(\varepsilon (T)\log T\rightarrow \infty $ as $%
T\rightarrow \infty )$ that 
\begin{equation*}
\Lambda \geq \sqrt{\frac{112}{12\left( 0.34613\right) \pi ^{2}}}=1.652\,9,
\end{equation*}%
which the desired value (\ref{h1}). The proof is complete.

\bigskip

Next in the following, we will derive an inequality from the Opial
inequality (\ref{Boyd}) which allows use to use the moments (\ref{hall2})
and (\ref{hall3}) ((\ref{Hu})) to derive a new unconditional (conditional)
lower bound for $\Lambda .$ As a special case of ((\ref{Boyd}) if $q=0,$ and 
$p=r=2k$, then the inequality (\ref{Boyd}) reduces to 
\begin{equation}
\int_{a}^{b}\left\vert y(t)\right\vert ^{2k}dt\leq A^{2k}(k)\left(
b-a\right) ^{2k}\int_{a}^{b}\left\vert y^{^{\prime }}(t)\right\vert ^{2k}dt,
\label{B1}
\end{equation}%
where 
\begin{equation}
A(k):=\left( \frac{2(k)}{2(k)-1}\right) ^{\frac{1}{2k}}\left( 2(k)\right) ^{%
\frac{2(k)-1}{2(k)}}\left( \Gamma (\frac{1}{2(k)})\Gamma (\frac{2(k)-1}{2(k)}%
)\right) ^{-1},  \label{AB}
\end{equation}%
and $\Gamma (u)$ is the Euler gamma function. The inequality (\ref{B1}) has
immediate application to the case where $y(a)=y(b)=0$. Choose $c=(a+b)/2$
and apply (\ref{B1}) to $[a,$ $c]$ and $[c,$ $b]$ and then add to obtain%
\begin{eqnarray}
\int_{a}^{b}\left\vert y(t)\right\vert ^{2k}dt &\leq &A^{2k}(k)\left( \frac{%
b-a}{2}\right) ^{2k}\left\{ \left( \int_{a}^{c}\left\vert y^{^{\prime
}}(t)\right\vert ^{2k}dt\right) +\left( \int_{c}^{b}\left\vert y^{^{\prime
}}(t)\right\vert ^{2k}dt\right) \right\}  \notag \\
&\leq &A^{2k}(k)\left( \frac{b-a}{2}\right) ^{2k}\int_{a}^{b}\left\vert
y^{^{\prime }}(t)\right\vert ^{2k}dt.  \label{B2.6}
\end{eqnarray}%
From this, we have 
\begin{equation*}
\int_{0}^{\pi }\left\vert y(t)\right\vert ^{2k}dt\leq A^{2k}(k)\left( \frac{%
\pi }{2}\right) ^{2k}\int_{0}^{\pi }\left\vert y^{^{\prime }}(t)\right\vert
^{2k}dt.
\end{equation*}%
with $y(0)=0=y(\pi ).$ From this inequality, we deduce that if $y\in
C^{1}[a,b]$\textit{\ }with $y(a)=0=y(b),$ then we have%
\begin{equation}
\int_{a}^{b}\left( \frac{b-a}{\pi }\right) ^{2k}\left( y^{^{\prime
}}(t)\right) ^{2k}dt\geq \left( \frac{2}{\pi }\right) ^{2k}\frac{1}{A^{2k}(k)%
}\int_{a}^{b}(y(t))^{2k}dt.  \label{B2}
\end{equation}%
Using the formula (\ref{AB}), we have the following values of $A(k)$ for $%
k=1,2,....15.$%
\begin{equation*}
\begin{tabular}[t]{|l|l|l|l|l|}
\hline
$A(1)$ & $A(2)$ & $A(3)$ & $A(4)$ & $A(5)$ \\ \hline
$0.63662$ & $0.68409$ & $0.73026$ & $0.76409$ & $0.7896$ \\ \hline
$A(6)$ & $A(7)$ & $A(8)$ & $A(9)$ & $A(10)$ \\ \hline
$0.80955$ & $0.825\,62$ & $0.83888$ & $0.85003$ & $0.85956$ \\ \hline
$A(11)$ & $A(12)$ & $A(13)$ & $A(14)$ & $A(15)$ \\ \hline
$0.8678$ & $0.87502$ & $0.88141$ & $0.88709$ & $0.89219$ \\ \hline
\end{tabular}%
\end{equation*}

\begin{center}
Table 2. The values of $A(k)$ for $k=1,2,...,15.$
\end{center}

In the following theorem, assuming the Riemann hypothesis, we apply the
inequality (\ref{B2}) when $k=2$ and using the moments (\ref{hall2}) and (%
\ref{hall3}) to derive an unconditional value for $\Lambda .$

\bigskip

\textbf{Theorem 2.2}.\textit{\ Let }$\varepsilon (T)\rightarrow 0$\textit{\
in such a way that }$\varepsilon (T)\log T\rightarrow \infty .$\textit{\
Then for sufficiently large }$T$\textit{, there exists an interval contained
in }$[T,(1+\varepsilon (T))T]$\textit{\ which is free of zeros of }$Z(t)$%
\textit{\ and having length at least }%
\begin{equation*}
\frac{1}{\pi \left( 0.68409\right) }\sqrt[4]{560}\left\{ 1+O\left( \frac{1}{%
\varepsilon (T)\log T}\right) \right\} \frac{2\pi }{\log T}.
\end{equation*}%
\textit{Thus }%
\begin{equation}
\Lambda \geq \frac{1}{\pi \left( 0.68409\right) }\sqrt[4]{560}=2.2635\text{.}
\label{BU}
\end{equation}%
\textbf{Proof. }As in the proof of Theorem 2.1,\textbf{\ }we follow the
arguments in \cite{Hall2} to prove our theorem. Suppose that $t_{l}$ is the
first zero of $Z(t)$ not less than $T$ and $t_{m}$ the last zero not greater
than $(1+\varepsilon )T$ where $\varepsilon (T)\rightarrow 0$ in such a way
that $\varepsilon (T)\log T\rightarrow \infty .$ Suppose further that for $%
l\leq n<m$, we have 
\begin{equation}
L_{n}=t_{n+1}-t_{n}\leq \frac{2\pi \kappa }{\log T}.  \label{B0}
\end{equation}%
Applying the inequality (\ref{B2}) with $y(t)=Z(t)$ and $k=2$, we have 
\begin{equation*}
\int_{t_{n}}^{t_{n+1}}\left[ \left( \frac{L_{n}}{\pi }\right) ^{4}\left(
Z^{^{\prime }}(t)\right) ^{4}-\left( \frac{2}{\pi A(2)}\right) ^{4}(Z(t))^{4}%
\right] dt\geq 0.
\end{equation*}%
Since the inequality remains true if we replace $L_{n}/\pi $ by $2\kappa
/\log T$, we have 
\begin{equation}
\int_{t_{n}}^{t_{n+1}}\left[ \left( \frac{2\kappa }{\log T}\right)
^{4}\left( Z^{^{\prime }}(t)\right) ^{4}-\left( \frac{2}{\pi A(2)}\right)
^{4}(Z(t))^{4}\right] dt\geq 0.  \label{B11}
\end{equation}%
Summing (\ref{B11}) over $n,$ using (\ref{hall2}) and (\ref{hall3}), we
obtain 
\begin{eqnarray*}
&&\frac{1}{1120\pi ^{2}}\left( \frac{2\kappa }{\log T}\right) ^{4}T\log
^{8}(T)+O(T\log ^{7}) \\
&&-\left( \frac{2}{\pi A(2)}\right) ^{4}\frac{1}{2\pi ^{2}}(T\log
^{4}T)+O(T\log ^{3}T) \\
&=&\frac{\left( 2\kappa \right) ^{4}}{1120\pi ^{2}}T\log ^{4}(T)+O(T\log
^{7}T) \\
&&-\left( \frac{2}{\pi A(2)}\right) ^{4}\frac{1}{2\pi ^{2}}(T\log
^{4}T)+O(T\log ^{3}T).
\end{eqnarray*}%
Follows the proof of Theorem 1 in \cite{Hall2}, we obtain 
\begin{equation*}
\kappa ^{4}\geq \frac{1}{2^{4}}\left( \frac{2}{\pi A(2)}\right) ^{4}\frac{%
1120\pi ^{2}}{2\pi ^{2}}+O(1/\varepsilon (T)\log T).
\end{equation*}%
Now, using the value of $A(2)$ from Table 2, we have (noting $(\varepsilon
(T)\log T\rightarrow \infty $ as $T\rightarrow \infty )$ that 
\begin{equation*}
\Lambda \geq \frac{1}{\pi }\frac{1}{\left( 0.68409\right) }\sqrt[4]{\frac{%
1120}{2}}=2.2635,
\end{equation*}%
which the desired value (\ref{BU}). The proof is complete.

\bigskip

In the following, we will establish some explicit formulae for the gaps
between the zeros of the Riemann zeta function and use them to find new
conditional series of lower bounds. First, we will apply the inequality (\ref%
{Kv3}). As usual, we assume that the Riemann hypothesis is true. As a
special case of (\ref{Kv3}), if $y(0)=y(\pi )=0$, $p=2k-2h$ and $q=2h$, we
have%
\begin{equation}
\int_{0}^{\pi }\left\vert y(t)\right\vert ^{2k-2h}\left\vert y^{^{\prime
}}(t)\right\vert ^{2h}dt\leq K(h,k)(\frac{\pi }{2})^{2k-2h}\left(
\int_{0}^{\pi }\left\vert y^{^{\prime }}(t)\right\vert ^{2k}dt\right) ,
\label{hk}
\end{equation}%
where 
\begin{equation}
K(h,k)=\frac{hk^{2k-1}}{\left( (k-h)L(k,h)+h\right) ^{2k}},\text{ \ }h,\text{
}k\neq 0\text{,}  \label{lk}
\end{equation}%
and 
\begin{equation}
L(h,k)=\int_{0}^{1}\left( \frac{1}{1-\lambda s^{2k-2h}}\right) ds\text{, \ \ 
}\lambda =\frac{k(2h-1)}{h(2k-1)}.  \label{lk1}
\end{equation}%
\textbf{\ Theorem 2.3. }\textit{On the hypothesis that the Riemann
hypothesis is true and (\ref{Hu}) is correctly predicted, we have\ }%
\begin{equation}
\Lambda \geq \Lambda ^{\ast }(h,k):=\frac{1}{\pi }\left( \frac{1}{K(h,k)}%
\frac{b(h,k)}{b(k,k)}\right) ^{\frac{1}{2k-2h}},\text{ }h\neq k\neq 0
\label{kli}
\end{equation}%
\textit{where }%
\begin{equation}
b(h,k):=\frac{2h!}{8^{h}h!}H(h,k)\left( \dprod_{j=0}^{k-1}\frac{j!}{(j+k)!}%
\right) .  \label{hk2}
\end{equation}%
\textbf{Proof.} As in the proof of Theorem 2.2 by\textbf{\ }applying the
inequality (\ref{hk}) with $y=Z(t)$, we have 
\begin{equation*}
\int_{t_{n}}^{t_{n+1}}\left[ \left( \frac{L_{n}}{\pi }\right) ^{2k}\left(
Z^{^{\prime }}(t)\right) ^{2k}-\frac{2^{2k-2h}}{K(h,k)\pi ^{2k-2h}}\left( 
\frac{L_{n}}{\pi }\right) ^{2h}\left\vert Z(t)\right\vert ^{2k-2h}\left\vert
Z^{^{\prime }}(t)\right\vert ^{2h}\right] dt\geq 0.
\end{equation*}%
Since the inequality remains true if we replace $L_{n}/\pi $ by $2\kappa
/\log T$, we have%
\begin{equation*}
\int_{t_{n}}^{t_{n+1}}\left( \frac{2\kappa }{\log T}\right) ^{2k}\left\vert
Z^{^{\prime }}(t)\right\vert ^{2k}
\end{equation*}%
\begin{equation*}
-\int_{t_{n}}^{t_{n+1}}\frac{2^{2k-2h}}{K(h,k)\pi ^{2k-2h}}\left( \frac{%
2\kappa }{\log T}\right) ^{2h}\left\vert Z(t)\right\vert ^{2k-2h}\left\vert
Z^{^{\prime }}(t)\right\vert ^{2h}dt\geq 0.
\end{equation*}%
Summing (\ref{BA3}) over $n$ and using (\ref{Hu}), we obtain%
\begin{equation*}
\left( \frac{2\kappa }{\log T}\right) ^{2k}a(k)b(k,k)T\left( \log T\right)
^{k^{2}+2k}
\end{equation*}%
\begin{equation*}
-\frac{2^{2k-2h}}{K(h,k)\pi ^{2k-2h}}\left( \frac{2\kappa }{\log T}\right)
^{2h}a(k)b(h,k)T\left( \log T\right) ^{k^{2}+2h}dt\geq 0.
\end{equation*}%
This implies that 
\begin{eqnarray*}
&&T\left( \log T\right) ^{k^{2}}\left\{ \left( 2\kappa \right)
^{2k}a(k)b(k,k)-\frac{2^{2k-2h}}{K(h,k)\pi ^{2k-2h}}\left( 2\kappa \right)
^{2h}a(k)b(h,k)\right\} \\
&\geq &o(T\left( \log T\right) ^{k^{2}}),
\end{eqnarray*}%
whence%
\begin{equation*}
\kappa ^{2k-2h}\geq \frac{1}{\pi ^{2k-2h}K(h,k)}\frac{b(h,k)}{b(k,k)}+o(1),%
\text{ (as }T\rightarrow \infty ).
\end{equation*}%
This implies that 
\begin{equation*}
\Lambda ^{2k-2h}(k)\geq \frac{1}{\pi ^{2k-2h}K(h,k)}\frac{b(h,k)}{b(k,k)},%
\text{ }h\neq k\neq 0.
\end{equation*}%
which is the desired inequality and completes the proof.

\bigskip

To apply the inequality (\ref{kli}), we will need the following values of $%
b(1,k)$ and $b(k,k)$ which are determined from (\ref{hk2}) where $H(h,k)$
are defined as in Table 1$:$

\begin{eqnarray*}
b(1,2) &=&\tfrac{1}{720},\text{ }b(2,2)=\tfrac{1}{6720}\text{, }b(1,3)=%
\tfrac{1}{1209\,600},\text{ }b(3,3)=\tfrac{1}{496742400}, \\
b(1,4) &=&\tfrac{1}{219469\,824\,000}\text{, \ \ }b(4,4)=\tfrac{31}{%
271159356948480000}, \\
b(1,5) &=&\tfrac{1}{8760533070643200\,000},\text{ }b(5,5)=\tfrac{227}{%
12854317559387145633792000000}, \\
b(1,6) &=&\tfrac{1}{127288\,050\,516627\,176\,816640\,000\,000}, \\
b(6,6) &=&\tfrac{133933}{25516459094444104187401241999966208000000000}, \\
b(1,7) &=&\tfrac{1}{998707926079695101611943783301120000000000}, \\
b(7,7) &=&\tfrac{2006509}{895370835179\,281010%
\,419215815294340559070476369920000000000000}.
\end{eqnarray*}

\begin{center}
Table 3: The values of $b(1,k)$ and $b(k,k)$ for $k=2,...,7.$
\end{center}

\bigskip

Also, we need the following values of $K(1,k)$ which are determined from the
formula (\ref{lk}) by using (\ref{lk1}) for $k=2,3,...,7:$%
\begin{eqnarray*}
K(1,2) &=&0.239\,61,\text{ }K(1,3)=0.161\,87,\text{ }K(1,4)=0.122\,27, \\
K(1,5) &=&9.823\,8\times 10^{-2},\text{ }K(1,6)=8.212\,8\times 10^{-2},\text{
}K(1,7)=0.07055.
\end{eqnarray*}%
Now, we are ready to derive a series of the lower bounds of $\Lambda (k)$
for $k=2,3,...,7.$ These lower bounds are determined by using the formula (%
\ref{kli}) and presented in the following table: 
\begin{equation*}
\begin{tabular}[t]{|l|l|l|l|l|l|}
\hline
$\Lambda (2)$ & $\Lambda (3)$ & $\Lambda (4)$ & $\Lambda (5)$ & $\Lambda (6)$
& $\Lambda (7)$ \\ \hline
$1.9866$ & $2.2591$ & $2.6407$ & $3.0208$ & $3.3800$ & $3.7124$ \\ \hline
\end{tabular}%
\end{equation*}

\begin{center}
Table 3: The lower bounds for $\Lambda (k)$ for $k=2,3,...,7$ by using the
formula (\ref{kli}).
\end{center}

In the following, we will apply the inequality (\ref{B2}) to establish a new
explicit formula for $\Lambda (k)$. As usual, we assume that the Riemann
hypothesis is true and the moments in (\ref{A1}) and (\ref{A2}) are
correctly predicted.

\bigskip

\textbf{Theorem 2.4}. \textit{Assuming the Riemann hypothesis and the
moments in (\ref{A1}) and (\ref{A2}) are correctly predicted, we have }%
\begin{equation}
\Lambda (k)\geq \frac{1}{\pi A(k)}\left( \frac{b_{k}}{c_{k}}\right) ^{\frac{1%
}{2k}},\text{ \ for }k\geq 3,  \label{BA1}
\end{equation}%
\textit{where }$A(k)$\textit{\ is defined as in (\ref{AB}).}

\textbf{Proof. }As in the proof of Theorem 2.2 by\textbf{\ }applying the
inequality (\ref{B2}) with $y=Z(t)$, we have 
\begin{equation*}
\int_{t_{n}}^{t_{n+1}}\left[ \left( \frac{L_{n}}{\pi }\right) ^{2k}\left(
Z^{^{\prime }}(t)\right) ^{2k}-\left( \frac{2}{\pi A(k)}\right)
^{2k}(Z(t))^{2k}\right] dt\geq 0.
\end{equation*}%
Since the inequality remains true if we replace $L_{n}/\pi $ by $2\kappa
/\log T$, we have 
\begin{equation}
\int_{t_{n}}^{t_{n+1}}\left[ \left( \frac{2\kappa }{\log T}\right)
^{2k}\left( Z^{^{\prime }}(t)\right) ^{2k}-\left( \frac{2}{\pi A(k)}\right)
^{2k}(Z(t))^{2k}\right] dt\geq 0.  \label{BA3}
\end{equation}%
Summing (\ref{BA3}) over $n,$ using (\ref{A1}) and (\ref{A2}), we obtain 
\begin{eqnarray*}
&&a_{k}c_{k}\left( \frac{2\kappa }{\log T}\right) ^{2k}T\left( \log T\right)
^{k^{2}+2k}-a_{k}b_{k}\left( \frac{2}{\pi A(k)}\right) ^{2k}T\left( \log
T\right) ^{k^{2}} \\
&=&\left( a_{k}c_{k}\kappa ^{2k}(2^{2k})-a_{k}b_{k}\left( \frac{2}{\pi A(k)}%
\right) ^{2k}\right) T\left( \log T\right) ^{k^{2}}\geq O(T\log ^{k^{2}}T),
\end{eqnarray*}%
whence%
\begin{equation*}
\kappa ^{2k}\geq \frac{a_{k}b_{k}}{2^{2k}a_{k}c_{k}}\frac{1}{A^{2k}(k)}=%
\frac{b_{k}}{2^{2k}c_{k}}\left( \frac{2}{\pi A(k)}\right) ^{2k},\text{ (as }%
T\rightarrow \infty ).
\end{equation*}%
This implies that 
\begin{equation*}
\Lambda ^{2k}(k)\geq \frac{b_{k}}{2^{2k}c_{k}}\left( \frac{2}{\pi A(k)}%
\right) ^{2k},
\end{equation*}%
and then we obtain the desired inequality (\ref{BA1}). The proof is complete.

\bigskip

Conrey, Rubinstein and Snaith \cite{C} gave some explicit values of the
parameter $c(k)$ for $k=1,2,...,15.$ Using the values of $c(k)$ due to
Conrey, Rubinstein and Snaith \cite{C} and the relation (\ref{B}), we have
the following values of $b(k)/c(k)$ for $k=1,2,...,15$ that will be used in
this paper:

\begin{center}
\begin{eqnarray*}
\frac{b_{1}}{c_{1}} &=&2^{2}\cdot 3,\text{ \ \ }\frac{b_{2}}{c_{2}}=\tfrac{%
2^{6}\cdot 3\cdot 5\cdot 7}{2^{2}3}\text{, \ }\frac{b_{3}}{c_{3}}=\tfrac{%
2^{12}\cdot 3^{2}\cdot 5^{2}\cdot 7^{2}\cdot 11}{2^{6}3^{3}5},\text{\ }\frac{%
b_{4}}{c_{4}}=\tfrac{2^{20}\cdot 3^{10}\cdot 5^{4}\cdot 7^{2}\cdot 11\cdot 13%
}{31\cdot 2^{12}3^{5}5^{3}7}, \\
\frac{b_{5}}{c_{5}} &=&\tfrac{2^{30}\cdot 3^{12}\cdot 5^{6}\cdot 7^{4}\cdot
11\cdot 13^{2}\cdot 17\cdot 19}{227\cdot 2^{20}3^{9}5^{5}7^{3}},\text{ }%
\frac{b_{6}}{c_{6}}=\tfrac{2^{42}\cdot 3^{19}\cdot 5^{9}\cdot 7^{6}\cdot
11^{3}\cdot 13^{3}\cdot 17\cdot 19\cdot 23}{67\cdot 1999\cdot
2^{30}3^{15}5^{7}7^{5}11}\text{,} \\
\frac{b_{7}}{c_{7}} &=&\tfrac{2^{56}\cdot 3^{28}\cdot 5^{13}\cdot 7^{8}\cdot
11^{4}\cdot 13^{3}\cdot 17^{2}\cdot 19^{2}\cdot 23}{43\cdot 46663\cdot
2^{42}3^{21}5^{9}7^{7}11^{3}13},\text{ }\frac{b_{8}}{c_{8}}=\tfrac{%
2^{72}\cdot 3^{34}\cdot 5^{16}\cdot 7^{11}\cdot 11^{6}\cdot 13^{4}\cdot
17^{3}\cdot 19^{2}\cdot 23\cdot 29\cdot 31}{46743947\cdot
2^{56}3^{28}5^{12}7^{9}11^{5}13^{3}}, \\
\frac{b_{9}}{c_{9}} &=&\tfrac{2^{90}\cdot 3^{42}\cdot 5^{21}\cdot
7^{14}\cdot 11^{8}\cdot 13^{6}\cdot 17^{3}\cdot 19^{3}\cdot 23^{2}\cdot
29\cdot 31}{19583\cdot 16249\cdot 2^{72}3^{36}5^{16}7^{11}11^{7}13^{5}17},%
\frac{b_{10}}{c_{10}}=\tfrac{2^{110}\cdot 3^{55}\cdot 5^{25}\cdot
7^{17}\cdot 11^{10}\cdot 13^{8}\cdot 17^{5}\cdot 19^{4}\cdot 23^{3}\cdot
29\cdot 31\cdot 37}{3156627824489\cdot
2^{90}3^{44}5^{20}7^{13}11^{9}13^{7}17^{3}19}, \\
\frac{b_{11}}{c_{11}} &=&\tfrac{2^{132}\cdot 3^{63}\cdot 5^{31}\cdot
7^{18}\cdot 11^{12}\cdot 13^{10}\cdot 17^{5}\cdot 19^{5}\cdot 23^{4}\cdot
29^{2}\cdot 31^{2}\cdot 37\cdot 41\cdot 43}{59\cdot 11332613\cdot 33391\cdot
2^{110}3^{53}5^{24}7^{16}11^{11}13^{9}17^{5}19^{3}},\text{ } \\
\frac{b_{12}}{c_{12}} &=&\tfrac{2^{156}\cdot 3^{75}\cdot 5^{37}\cdot
7^{23}\cdot 11^{15}\cdot 13^{12}\cdot 17^{8}\cdot 19^{7}\cdot 23^{4}\cdot
29^{3}\cdot 31^{2}\cdot 41\cdot 43\cdot 47}{241\cdot 251799899121593\cdot
2^{132}3^{63}5^{28}7^{20}11^{13}13^{11}17^{7}19^{5}23}, \\
\frac{b_{13}}{c_{13}} &=&\tfrac{2^{182}\cdot 3^{90}\cdot 5^{42}\cdot
7^{28}\cdot 11^{17}\cdot 13^{14}\cdot 17^{10}\cdot 19^{8}\cdot 23^{5}\cdot
29^{3}\cdot 31^{3}\cdot 37^{2}\cdot 41\cdot 43\cdot 47}{285533\cdot
37408704134429\cdot 2^{156}3^{73}5^{34}7^{24}11^{15}13^{13}17^{9}19^{7}23^{3}%
}\text{, } \\
\frac{b_{14}}{c_{14}} &=&\tfrac{2^{210}\cdot 3^{100}\cdot 5^{50}\cdot
7^{31}\cdot 11^{20}\cdot 13^{17}\cdot 17^{12}\cdot 19^{10}\cdot 23^{7}\cdot
29^{4}\cdot 31^{4}\cdot 37^{2}\cdot 41\cdot 43\cdot 47\cdot 53}{197\cdot
1462253323\cdot 6616773091\cdot
2^{182}3^{86}5^{42}7^{28}11^{17}13^{15}17^{11}19^{9}23^{5}}, \\
\frac{b_{15}}{c_{15}} &=&\tfrac{2^{240}\cdot 3^{117}\cdot 5^{57}\cdot
7^{37}\cdot 11^{22}\cdot 13^{19}\cdot 17^{14}\cdot 19^{11}\cdot 23^{9}\cdot
29^{5}\cdot 31^{5}\cdot 37^{3}\cdot 41^{2}\cdot 43^{2}\cdot 47\cdot 53\cdot
59}{1625537582517468726519545837\cdot
2^{210}3^{102}5^{50}7^{32}11^{19}13^{17}17^{13}19^{11}23^{7}29}.
\end{eqnarray*}

Table 4. The values of $b(k)/c(k)$ for $k=1,2,...,15.$
\end{center}

Using the formula (\ref{BA1}), the values of $b(k)/c(k)$ in Table 4, and the
values of $A(k)$ in Table 2, we have the following new lower bounds for $%
\Lambda (k)$ for $k=1,2,...,15$:

\begin{equation*}
\begin{tabular}[t]{|l|l|l|l|l|}
\hline
$\Lambda (1)$ & $\Lambda (2)$ & $\Lambda (3)$ & $\Lambda (4)$ & $\Lambda (5)$
\\ \hline
$1.732\,1$ & $2.2635$ & $2.7080$ & $3.1257$ & $3.5177$ \\ \hline
$\Lambda (6)$ & $\Lambda (7)$ & $\Lambda (8)$ & $\Lambda (9)$ & $\Lambda
(10) $ \\ \hline
$3.8813$ & $4.2150$ & $4.5196$ & $4.7985$ & $5.0560$ \\ \hline
$\Lambda (11)$ & $\Lambda (12)$ & $\Lambda (13)$ & $\Lambda (14)$ & $\Lambda
(15)$ \\ \hline
$5.2962$ & $5.5225$ & $5.7373$ & $5.9424$ & $6.1392$ \\ \hline
\end{tabular}%
\end{equation*}

\begin{center}
Table 4. The lower bounds for $\Lambda (k)$ for $k=1,2,...,15.$
\end{center}

From this table we have the following theorem.

\textbf{Theorem 2.5. }\textit{On the hypothesis that the Riemann hypothesis
is true, (\ref{A1}) and (\ref{A2}) are correctly predicted, we have\ }%
\begin{equation}
\Lambda \geq 6.1392  \label{BA4}
\end{equation}

\begin{remark}
Using the explicit formulae for the $b(k)$ and $c(k)$ (which would via (\ref%
{BA1}) help to decide whether the conjecture $\Lambda =\infty $ is true
subject to the Riemann hypothesis), we have the following formula%
\begin{eqnarray*}
\Lambda (k) &\geq &\frac{1}{\pi A(k)}\left( \dprod_{j=0}^{k-1}\frac{j!}{%
(j+k)!}\right) ^{\frac{1}{2k}} \\
&&\times (-1)^{\frac{k(k+1)}{2}}\dsum_{m\in P_{O}^{k+1}(2k)}^{k}\left( 
\begin{array}{c}
2k \\ 
m%
\end{array}%
\right) \left( \frac{-1}{2}\right) ^{m_{0}}\left( \dprod_{i=1}^{k}\tfrac{1}{%
(2k-i+m_{i})!}\right) M_{i,j},
\end{eqnarray*}%
where $\ A(k)$ is defined as in (\ref{AB}), $M_{i,j}$ as defined in (\ref{M}%
) and $P_{O}^{k+1}(2k)$ denotes the set of partitions $m=(m_{0},...,m_{k})$
of $2k$ into nonnegative parts.
\end{remark}

\begin{remark}
The lower bound in (\ref{BA4}) means that consecutive nontrivial zeros often
differ by at least $6.1392$ times the average spacing. This value improves
the value of $\Lambda \geq 4.71474396$ that has been obtained by Hall.
\end{remark}

\end{document}